\documentclass[fleqn]{mat01}
\usepackage{times,mathtimy,amssymb,latexsym}
\begin{document}

\setcounter{page}{13} \firstpage{13}

\newtheorem{theor}{\bf Theorem}
\newtheorem{coro}[theor]{\rm COROLLARY}
\newtheorem{propo}[theor]{\rm PROPOSITION}
\newtheorem{lem}[theor]{Lemma}
\newtheorem{rema}[theor]{Remark}
\newtheorem{exam}[theor]{Example}

\title{Khinchin's inequality, Dunford--Pettis and compact operators
on the space $\pmb{C([0,1],X)}$}

\markboth{Dumitru Popa}{Banach spaces of continuous functions}

\author{DUMITRU POPA}

\address{Department of Mathematics, University of Constanta,
8700~Constanta, Romania\\
\noindent E-mail: dpopa@univ-ovidius.ro}

\volume{117}

\mon{February}

\parts{1}

\pubyear{2007}

\Date{MS received 10 November 2005}

\begin{abstract}
We prove that if $X,Y$ are Banach spaces, $\Omega$ a compact
Hausdorff space and $U\hbox{\rm :}\ C(\Omega,X)\rightarrow Y$ is a
bounded linear operator, and if $U$ is a Dunford--Pettis operator
the range of the representing measure $G(\Sigma) \subseteq
DP(X,Y)$ is an uniformly Dunford--Pettis family of operators and
$\|G\|$ is continuous at $\emptyset$. As applications of this
result we give necessary and/or sufficient conditions that some
bounded linear operators on the space $C([0,1],X)$ with values in
$c_{0}$ or $l_{p}$, ($1\leq p<\infty$) be Dunford--Pettis and/or
compact operators, in which, Khinchin's inequality plays an
important role.
\end{abstract}


\maketitle

\noindent Let $\Omega$ be a compact Hausdorff space, $X$ a Banach
space, and $C(\Omega,X)$ the Banach space of continuous $X$-valued
functions on $\Omega$ under the uniform norm and $C(\Omega)$ when
$X$ is the scalar field. It is well-known (see \S1 of \cite{1},
Theorem~2.2 of \cite{3} (Representation Theorem) or Theorem
(Dinculeanu-Singer), p.~182 of \cite{6}) that if $Y$ is a Banach
space then any bounded linear operator $U\hbox{\rm :}\ C(\Omega
,X)\rightarrow Y$ has a finitely additive vector measure
$G\hbox{\rm :}\ \Sigma \rightarrow L(X,Y^{**})$, where $\Sigma$ is
the $\sigma $-field of Borel subsets of $\Omega$, such that
$y^{*}U(f) = \int_{\Omega}f\hbox{d}G_{y^{*}}$, $f\in C(\Omega
,X)$, $y^{*}\in Y^{*}$. The measure $G$ is called the representing
measure of $U$.

Also, for a bounded linear operator $U\hbox{\rm :}\
C(\Omega,X)\rightarrow Y$ we can associate in a natural way two
bounded linear operators $U^{\#}\hbox{\rm :}\ C(\Omega)\rightarrow
L(X,Y)$ and $U_{\#}\hbox{\rm :}\ X$ $\rightarrow L(C(\Omega),Y)$
defined by $(U^{\#}\varphi) (x) =U(\varphi \otimes x)$ and
$(U_{\#}x) (\varphi) =U(\varphi \otimes x)$, where for $\varphi
\in C(\Omega)$, $x\in X$ we denote $(\varphi \otimes x) (\omega)
=\varphi (\omega) x$.

For a $\sigma$-algebra $\Sigma \subseteq \mathcal{P}(S)$, $X$ a
Banach space and a vector measure $G\hbox{\rm :}\ \Sigma
\rightarrow X$, we denote $\tilde{G}(E) =\sup \{\|G(A)\||A\in
\Sigma, A\subseteq E\}$ the quasivariation of $G$, by $|G|$ and
$\|G\|$ the variation and semivariation of $G$ and we use the fact
that $\tilde{G}(E)\leq \|G\|(E) \leq 4\tilde{G}(E)$ for any $E\in
\Sigma$ (see chapter~I, Proposition~11, p.~4 of \cite{6}).

We denote by $B(\Sigma,X)$ the space of all totally measurable
functions endowed with the supnorm.

Also, for $[0,1]$ we denote by $\Sigma $ the $\sigma $-field of
Borel subsets, $\mu\hbox{\rm :}\ \Sigma \rightarrow [0,1]$ is the
Lebesgue measure and $(r_{n})_{n\in \Bbb{N}}$ is the sequence of
Rademacher functions.

If $\nu \in rcabv(\Sigma)$, $f\hbox{\rm :}\ [0,1] \rightarrow
\Bbb{K}$ is $\nu$-integrable and $\alpha\hbox{\rm :}\ \Sigma
\rightarrow \Bbb{K}$ is defined by $\alpha (E) =\int_{E}f(t)
\hbox{d}\nu (t)$, then $\int_{0}^{1}|f(t)|\hbox{d}|\nu|(t)
=|\alpha|([0,1])$ and $\tilde{\alpha}([0,1]) \leq
\int_{0}^{1}|f(t)|\hbox{d}|\nu|(t) \leq 4 \tilde{\alpha}([0,1])$.

If $G\hbox{\rm :}\ \Sigma \rightarrow L(X,Y) $ is a vector measure
we denote by $\|G\|$ the semivariation of $G$ defined by $\|G\|(E)
=\sup \{|G_{y^{*}}|(E)|\|y^{*}\|\leq 1\}$, $E\in \Sigma$, where
$G_{y^{*}}(E)=\langle G(E)x,y^{*}\rangle$ and we say that the
semivariation $\|G\|$ is continuous at $\emptyset $ if
$\|G\|(E_{k}) \rightarrow 0$ for $E_{k}\searrow \emptyset$,
$(E_{k})_{k\in \Bbb{N}}\subset \Sigma$. As is well-known, $\|G\|$
is continuous at $\emptyset$ if and only if there exists $\alpha
\geq 0$ a Borel measure on $\Sigma$ such that $\lim_{\alpha
(E)\rightarrow 0}\|G(E) \|=0$. Also, $|G|$ is the variation of $G$
and for $x\in X$ we write $G_{x}\hbox{\rm :}\ \Sigma \rightarrow
Y$ defined by $G_{x}(E) =G(E) (x)$ and if $\lambda\hbox{\rm :}\
\Sigma \rightarrow X^{*}$ is a vector measure for $x\in X$ we
write $\lambda x\hbox{\rm :}\ \Sigma \rightarrow \Bbb{K}$ defined
by $(\lambda x) (E) =\lambda (E) (x)$.

As is well-known, (see chapter~2, p.~32 of \cite{5}) if $X$ is a
Banach space, $1\leq p<\infty$ and $(x_{n}^{*})_{n\in
\Bbb{N}}\subseteq X^{*}$ is such that for any $x\in X$ the series
$\sum_{n=1}^{\infty}|x_{n}^{*}(x)|^{p}$ is convergent. Then
\begin{equation*}
w_{p}(x_{n}^{*}| n\in \Bbb{N}) =\sup\limits_{\|x\|\leq 1}
\left(\sum\limits_{n=1}^{\infty}|x_{n}^{*}(x)|^{p}
\right)^{\frac{1}{p}}<\infty
\end{equation*}
and we denote by $w_{p}(X^{*})$ the set of all such sequences.

Observe that if $p=1$, then $(x_{n}^{*})_{n\in \Bbb{N}}\in
w_{1}(X^{*}) $ if and only if $\sum_{n=1}^{\infty}x_{n}^{*}$ is a
weakly Cauchy series in $X^{*}$.

We recall that if $X$ and $Y$ are Banach spaces, a bounded linear
operator $U\hbox{\rm :}\ X\rightarrow Y$ is called Dunford--Pettis
operator if and only if for any $x_{n}\rightarrow 0$ {\it weak} it
follows that $U(x_{n}) \rightarrow 0$ in norm. We denote by
$DP(X,Y)$ the space of all the Dunford--Pettis operators from $X$
into $Y$. A~Banach space has the Schur property if the identity
operator is Dunford--Pettis.

We also need the following characterization of $weak$ convergence
in a $C(\Omega,X)$ space which will be used later without an
explicit reference. If $(f_{n})_{n\in \Bbb{N}}\subseteq
C(\Omega,X)$, then $f_{n}\rightarrow 0$ {\it weak} if and only if
$\sup_{n\in \Bbb{N},\omega \in \Omega} \|f_{n}(\omega)\|<\infty$
and $f_{n}(\omega) \rightarrow 0$ {\it weak} for any $\omega \in
\Omega$ (see Theorem~2 of \cite{2}).

For $c_{0}$ or $l_{p}$ with $1\leq p\leq \infty $ we denote by
$e_{n}$ the standard unit vectors in these spaces.

All notations and notions used and not defined in this paper are
either standard or can be found in \cite{5} or \cite{6}.

In Theorem~3.1 of \cite{10} it is proved that if $U\hbox{\rm :}\
C(\Omega,X) \rightarrow Y$ is a Dunford--Pettis operator, then the
representing measure has the property that $G(E) \in DP(X,Y)$ for
any $E\in \Sigma$ and $\|G\|$ is continuous at $\emptyset$ and
that this condition is necessary and sufficient if and only if $X$
has the Schur property.

We will prove in theorem~4 below, that Dunford--Pettis operators
on $C(\Omega,X) $ satisfies a much stronger condition. In order to
prove this result we introduce the following notion.

Let $X$ be a Banach space and $(Y_{i})_{i\in I}$ a family of
Banach spaces.

A family $\{ U_{i}\in L(X,Y_{i})| i\in I\} $ is said to be {\it
uniformly Dunford--Pettis family of operators} if and only if for
any $x_{n}\rightarrow 0$ {\it weak} it follows that $\sup_{i\in
I}\|U_{i}(x_{n}) \|\rightarrow 0$.

When there is some risk of confusion we write $\|\cdot \text{}
\|_{i}$ for the norm in the Banach space $Y_{i}$.

In the sequel we give necessary and sufficient conditions that a
sequence of bounded linear operators be an uniformly
Dunford--Pettis family.

\begin{propo}$\left.\right.$\vspace{.5pc}

\noindent Let $X$ be a Banach space{\rm ,} $(Y_{n})_{n\in
\Bbb{N}}$ a sequence of Banach spaces and $U_{n}\in L(X,Y_{n})$
for any $n\in \Bbb{N}$. The following assertions are
equivalent{\rm :}
\begin{enumerate}
\renewcommand\labelenumi{\rm (\roman{enumi})}
\leftskip .4pc
\item $(U_{n})_{n\in \Bbb{N}}$ is an uniformly Dunford--Pettis
family.

\item $U_{n}\in DP(X,Y_{n})$ for any $n\in \Bbb{N}$ and
for any sequence $x_{n}\rightarrow 0$ weak it follows that
$\|U_{n}(x_{n})\|_{n}\rightarrow 0$ in norm i.e. the diagonal
sequence of the matrix $(U_{n}(x_{k}))_{n,k\in \Bbb{N}}$ is null
convergent.

\item $\sup_{n\in \Bbb{N}}\|U_{n}\|_{n}<\infty$ and
the operator $U\hbox{\rm :}\ X\rightarrow l_{\infty}(Y_{n}| n\in
\Bbb{N})$ defined by $U(x) = (U_{n}(x))_{n\in \Bbb{N}}$ is
Dunford--Pettis.
\end{enumerate}
\end{propo}

\begin{proof}$\left.\right.$

\noindent (i) or (iii)$\Rightarrow$(ii).\ \ It is trivial.

\noindent (ii)$\Rightarrow$(i).\ \ Indeed, it is easy to see that
from (ii) it follows that for any sequence $x_{n}\rightarrow 0$
{weak} and any two subsequences $(k_{n})_{n\in \Bbb{N}}$ and
$(p_{n})_{n\in \Bbb{N}} $ of $\Bbb{N}$ it follows that
$\|U_{k_{n}}(x_{p_{n}})\|_{k_{n}}\rightarrow 0$. Since $U_{k}$ is
a Dunford--Pettis operator for any $k\in \Bbb{N}$, we have
$\lim_{n\rightarrow \infty}\|U_{k}(x_{n}) \|_{k}=0$. Now using the
well-known fact that if $(a_{nk})_{n,k\in \Bbb{N}}\subseteq
[0,\infty)$ is a double indexed sequence such that for any $k\in
\Bbb{N}$ we have $\lim_{n\rightarrow \infty}a_{kn}=0$, then
$\lim_{n\rightarrow \infty }a_{kn}=0$ uniformly in $k\in \Bbb{N}$
if and only if for any two subsequences $(k_{n})_{n\in \Bbb{N}}$
and $(p_{n})_{n\in \Bbb{N}}$ of $\Bbb{N}$ it follows that
$\lim_{n\rightarrow \infty}a_{k_{n}p_{n}}=0$, we deduce
$\sup_{k\in \Bbb{N}}\|U_{k}(x_{n}) \|_{k}\rightarrow 0$, i.e.
(ii).

\noindent (ii)$\Rightarrow$(iii). Let $c_{0}^{\rm weak}(X) =\{
(x_{n})_{n\in \Bbb{N}}\subseteq X| \text{}x_{n}\rightarrow 0\
\hbox{weak}\}$ which is a linear space for the natural operations
for addition and scalar multiplication and a Banach space for the
norm $\|(x_{n})_{n\in \Bbb{N}}\|=\sup_{n\in \Bbb{N}}\|x_{n}\|$.
Then (ii) affirms that the mapping $h\hbox{\rm :}\ c_{0}^{\rm
weak}(X) \rightarrow c_{0}(Y_{n}| n\in \Bbb{N}) $ defined by
$h((x_{n})_{n\in \Bbb{N}}) =(U_{n}(x_{n}))_{n\in \Bbb{N}}$ takes
its values in $c_{0}(Y_{n}| n\in \Bbb{N}) $ and, by an easy
application of the closed graph theorem, $h$ is bounded linear.
Then for any $n\in \Bbb{N}$ and $x\in X$ we have $\|U_{n}(x)
\|_{n}=\|h(0,\dots,0,x,0,\dots) \|\leq \|h\|\|x\|$ i.e. the family
$(U_{n})_{n\in \Bbb{N}}$ is pointwise bounded and thus uniformly
bounded, (by the uniform boundedness principle) i.e. $\sup_{n\in
\Bbb{N}}\|U_{n}\|_{n}<\infty$. Then the operator $U$ in (iii) is
well-defined, bounded linear and by the equivalence between (i)
and (ii), it follows that $U$ is\break Dunford--Pettis.
\end{proof}

As a consequence, from Proposition~1 we give a necessary and
sufficient condition that an operator with values in $c_{0}$ or
$l_{\infty}$ be Dunford--Pettis, completing a result from
Exercise~4, p.~114 of \cite{4}. Probably, this result is
well-known, but we do not know a\break reference.

\begin{coro}$\left.\right.$\vspace{.5pc}

\noindent Let $X$ be a Banach space{\rm ,} $(x_{n}^{*})_{n\in
\Bbb{N}} \subseteq X^{*}$ such that either $(x_{n}^{*})_{n\in
\Bbb{N}}$ is bounded or{\rm ,} $x_{n}^{*}\rightarrow 0\
\hbox{weak}^{*}$ and $U\hbox{\rm :}\ X\rightarrow l_{\infty}$ or
$c_{0}$ defined by $U(x) =(x_{n}^{*}(x))_{n\in \Bbb{N}}$.

Then the following assertions are equivalent{\rm :}
\begin{enumerate}
\renewcommand\labelenumi{\rm (\roman{enumi})}
\leftskip .2pc
\item $U$ is a Dunford--Pettis operator.

\item For any sequence $x_{n}\rightarrow 0$ weak it
follows that $x_{n}^{*}(x_{n}) \rightarrow 0$.\vspace{-.5pc}
\end{enumerate}
\end{coro}

In the next proposition the point (a)~is an extension of the
implication (ii) $\Rightarrow$ (i) in Theorem~3.1 in \cite{10}.

\begin{propo}$\left.\right.$\vspace{.5pc}

\noindent Let $X$ be a Banach space{\rm ,} $(Y_{n})_{n\in
\Bbb{N}}$ a sequence of Banach spaces and $T_{n},V_{n}\in
L(X,Y_{n})$ two sequences with $\sup_{n\in
\Bbb{N}}\|T_{n}\|_{n}<\infty$ and $\sup_{n\in
\Bbb{N}}\|V_{n}\|_{n}<\infty$. Let $U\hbox{\rm :}\ C([0,1],X)
\rightarrow c_{0}(Y_{n}|n\in \Bbb{N})$ be the operator defined by
\begin{equation*}
U(f) = \left( \int_{0}^{1} (T_{n}(f(t)) \sin 2\pi nt+ V_{n}(f(t))
\cos 2\pi nt) {\rm d}t\right)_{n\in \Bbb{N}}.
\end{equation*}
Then
\begin{enumerate}
\renewcommand\labelenumi{\rm (\alph{enumi})}
\leftskip .1pc
\item $U$ is Dunford--Pettis $\Leftrightarrow$ $(T_{n})_{n\in
\Bbb{N}}$ and $(V_{n})_{n\in \Bbb{N}}$ are uniformly
Dunford--Pettis.

\item $U$ is compact $\Leftrightarrow$ $T_{n}$
and $V_{n}$ are compact for any $n\in \Bbb{N}$ and
$\|T_{n}\|\rightarrow 0$ and $\|V_{n}\|\rightarrow 0$.
\end{enumerate}
\end{propo}

\begin{proof} The fact that $U$ takes its values in
$c_{0}(Y_{n}| n\in \Bbb{N}) $ follows from hypothesis, the
well-known fact that for any $f\in C[0,1]$ we have
$\int_{0}^{1}f(t) \sin 2\pi nt\hbox{d}t\rightarrow 0$ and
$\int_{0}^{1}f(t) \cos 2\pi nt\hbox{d}t\rightarrow 0$ and the
density of $C[0,1]\otimes X$ in $C([0,1],X)$. For any $n\in
\Bbb{N}$, let $S_{n}\hbox{\rm :}\ C([0,1],X) \rightarrow Y_{n}$ be
defined by
\begin{equation*}
S_{n}(f)=\int_{0}^{1}(T_{n}(f(t)) \sin 2\pi nt+V_{n}(f(t)) \cos
2\pi nt) \hbox{d}t.
\end{equation*}

(a) Suppose $U$ is a Dunford--Pettis operator. Let $k\in \Bbb{N}$
be fixed.

If $x_{n}\rightarrow 0$ weak, then $x_{n}\sin 2\pi kt\rightarrow
0$ weak in $C([0,1],X)$ and thus $U(x_{n}\sin 2\pi kt) \rightarrow
0$ in norm. If $U(x_{n}\sin 2\pi kt) =\frac{1}{2}
(0,\dots,0,T_{k}(x_{n}),0,\dots)$ then $T_{k}(x_{n}) \rightarrow
0$ in norm i.e. $T_{k}$ is Dunford--Pettis. Also, if
$x_{n}\rightarrow 0$ weak, then $x_{n}\sin 2\pi nt\rightarrow 0$
weak in $C([0,1],X)$, thus $U(x_{n}\sin 2\pi nt) \rightarrow 0$ in
norm and since $U(x_{n}\sin 2\pi nt)
=\frac{1}{2}(0,\dots,0,T_{n}(x_{n}),0,\dots)$ it follows that
$T_{n}(x_{n}) \rightarrow 0$ in norm. From Proposition~1 it
follows that $(T_{n})_{n\in \Bbb{N}}$ is an uniformly
Dunford--Pettis family. In the same way it can be proved that
$(V_{n})_{n\in \Bbb{N}}$ is an uniformly Dunford--Pettis family.

Suppose now that $(T_{n})_{n\in \Bbb{N}}$ and $(V_{n})_{n\in
\Bbb{N}}$ are uniformly Dunford--Pettis. Then, by the ideal
property of the class of all Dunford--Pettis operators, it follows
that $S_{n}$ is Dunford--Pettis for any $n\in \Bbb{N}$. Let
$(f_{n})_{n\in \Bbb{N}}\subseteq C([0,1],X) $ be such that
$f_{n}\rightarrow 0$ weak. Then $\sup_{n\in \Bbb{N} \text{,}t\in
[0,1]}\|f_{n}(t) \|<\infty $ and $f_{n}(t) \rightarrow 0$ weak for
any $t\in [0,1]$. By Proposition~1, $\|T_{n}(f_{n}(t)) \sin 2\pi
nt\|_{n}\rightarrow 0$ for any $t\in [0,1]$ and obviously
$\sup_{n\in \Bbb{N},t\in [0,1]} \|T_{n}(f_{n}(t)) \sin 2\pi
nt\|_{n}<\infty$. Now by Bartle's convergence theorem (p.~56 of
\cite{6}), it follows that $\int_{0}^{1}\|T_{n}(f_{n}(t)) \sin
2\pi nt\|_{n}\hbox{d}t\rightarrow 0$. Analogously
$\int_{0}^{1}\|V_{n}(f_{n}(t)) \cos 2\pi
nt\|_{n}\hbox{d}t\rightarrow 0$. By Proposition~1, $U$ is
Dunford--Pettis.

(b) Suppose $U$ is compact. Then (see Exercise~4, p.~144 of
\cite{4}) there is $0\leq \lambda_{n}\rightarrow 0$ such that
$\|S_{n}(f) \|\leq \lambda_{n}$ for any $f\in C([0,1],X) $ with
$\|f\|\leq 1$ and any $n\in \Bbb{N}$.

In particular, for any $x\in B_{X}$ and any $n\in \Bbb{N}$ we have
\begin{equation*}
\left\|\int_{0}^{1} (T_{n}(x) \sin^{2}2\pi nt+V_{n}(x) \sin 2\pi
nt\cos 2\pi nt) \hbox{d}t \right\| \leq \lambda_{n}
\end{equation*}
i.e. $\|T_{n}\|\leq 2\lambda_{n}$ for any $n\in \Bbb{N}$ and thus
$\|T_{n}\|\rightarrow 0$. Similarly $\|V_{n}\|\rightarrow 0$.

Also, by the ideal property of compact operators we obtain that
for any $n\in \Bbb{N}$ the operator $S_{n}$ is compact, in
particular the set
\begin{equation*}
\left\{ \int_{0}^{1}(T_{n}(x) \sin^{2}2\pi nt+V_{n}(x) \sin 2\pi
nt\cos 2\pi nt) \hbox{d}t| \|x\|\leq 1 \right\} \subseteq Y_{n}
\end{equation*}
is relatively norm compact, $T_{n}(B_{X}) $ is relatively norm
compact i.e. $T_{n}$ is compact. Analogously, $V_{n}$ is compact.

Conversely, by the ideal property of compact operators, it follows
that all $S_{n}$ are compact and also $\|S_{n}\|\leq \|T_{n}\|+$
$\|V_{n}\|\rightarrow 0$ i.e. $U$ is compact.
\end{proof}

The following theorem, which is the main result of our paper, is
an extension of Theorem~3.1 in \cite{10}.

\begin{theor}[\!]
Let $X,Y$ be Banach spaces{\rm ,} $\Omega$ a compact Hausdorff
space and $U\hbox{\rm :}\ C(\Omega,X)$ $\rightarrow Y$ a bounded
linear operator with $G$ its representing measure.

If $U$ is a Dunford--Pettis operator{\rm ,} then the range of the
representing measure $G(\Sigma) \subseteq DP(X,Y)$ is an uniformly
Dunford--Pettis family of operators and $\|G\|$ is continuous at
$\emptyset${\rm ,} or equivalently{\rm ,} for any
$x_{n}\rightarrow 0$ weak it follows that
$\widetilde{G_{x_{n}}}(\Omega) \rightarrow 0$ and $\|G\|$ is
continuous at $\emptyset$.
\end{theor}

\begin{proof}
If $U$ is a Dunford--Pettis operator, then clearly for any $x\in
X$ we have $U_{\#}(x) \in DP(C(\Omega),Y)$. Let $x_{n}\rightarrow
0$ weak. For $n\in \Bbb{N}$, let $\varphi _{n}\in C(\Omega)$ with
$\|\varphi_{n}\|\leq 1$ such that
\begin{equation*}
\|U_{\#}(x_{n})\| - \frac{1}{n}< \|U_{\#}(x_{n}) (\varphi_{n})\| =
\|U(\varphi_{n}\otimes x_{n})\|.
\end{equation*}

For any $\omega \in \Omega $ and any $x^{*}\in X^{*}$ we have
\begin{equation*}
|x^{*}(\varphi_{n}(\omega) x_{n})|\leq|x^{*}(x_{n})|\rightarrow 0
\ \ \hbox{and}\ \ \sup\limits_{n\in \Bbb{N}}\|\varphi_{n}\otimes
x_{n}\|\leq \sup\limits_{n\in \Bbb{N}}\|x_{n}\|<\infty,
\end{equation*}
hence $\varphi_{n}\otimes x_{n}\rightarrow 0$ weak. Since $U$ is
Dunford--Pettis we have $\|U(\varphi_{n}\otimes x_{n})
\|\rightarrow 0$ and thus $\|U_{\#}(x_{n}) \|\rightarrow 0$, which
means that $U_{\#}\hbox{\rm :}\ X\rightarrow DP(C(\Omega),Y) $ is
a Dunford--Pettis operator. Because for any $x\in X$ the operator
$U_{\#}(x)\hbox{\rm :}\ C[0,1] \rightarrow Y$ has the representing
measure $G_{x}\hbox{\rm :}\ \Sigma \rightarrow Y$ and
$\widetilde{G_{x}}(\Omega) \leq \|U_{\#}(x) \|=\|G_{x}\|(\Omega)
\leq 4\widetilde{G_{x}}(\Omega) $ we get that for any
$x_{n}\rightarrow 0$ weak it follows that $\sup_{E\in
\Sigma}\|G(E) (x_{n}) \|=\widetilde{G_{x_{n}}} (\Omega)
\rightarrow 0$.

The fact that $\|G\|$ is continuous at $\emptyset$ is proved in
Theorem~3.1 of \cite{10}.

We observe that the above proof is an obvious modification of the
proof of Proposition~7 in \cite{8}.
\end{proof}

Now we analyze the case of operators with values in $c_{0}$.

\begin{theor}[\!]$\left.\right.$
\begin{enumerate}
\renewcommand\labelenumi{\rm (\roman{enumi})}
\leftskip .2pc
\item Let $X$ be a Banach space{\rm ,} $(\lambda_{n})_{n\in
\Bbb{N}}\subseteq rcabv(\Sigma,X^{*})$ such that $\lambda_{n}(E)
\rightarrow 0$ $\hbox{weak}^{*}$ for any $E\in \Sigma$ and let
$U\hbox{\rm :}\ C([0,1],X) \rightarrow c_{0}$ be the operator
defined by
\begin{equation*}
\hskip -1.25pc U(f) =\left(\int_{0}^{1}f(t) {\rm d}\lambda_{n}(t)
\right)_{n\in \Bbb{N}}.
\end{equation*}

If $U$ is a Dunford--Pettis operator{\rm ,} then for any
$x_{n}\rightarrow 0$ weak it follows that $\sup_{k\in
\Bbb{N}}|\lambda_{k}x_{n}|([0,1]) \rightarrow 0$ and
$(\lambda_{n})_{n\in \Bbb{N}}$ is uniformly countably additive.

\item Let $X$ be a Banach space{\rm ,} $(\varphi_{n})_{n\in
\Bbb{N}}\subseteq B(\Sigma,X^{*})$ with $\sup_{n\in
\Bbb{N}\text{,}t\in [0,1]}\|\varphi_{n}(t) \|=M<\infty${\rm ,}
$(\nu_{n})_{n\in \Bbb{N}}\subseteq rcabv(\Sigma)$ an uniformly
countably additive pointwise bounded family such that for any
$E\in \Sigma$ and any $x\in X$ we have $\int_{E}\varphi_{n}(t) (x)
{\rm d}\nu_{n}(t) \rightarrow 0$. Let $U\hbox{\rm :}\ C([0,1],X)
\rightarrow c_{0}$ be the operator defined by
\begin{equation*}
\hskip -1.25pc U(f) =\left(\int_{0}^{1}\varphi_{n}(t) (f(t)) {\rm
d}\nu_{n}(t) \right)_{n\in \Bbb{N}}.
\end{equation*}

\leftskip .4pc (a) If $U$ is a Dunford--Pettis operator{\rm ,}
then for any $x_{n}\rightarrow 0$ weak it follows that
$\int_{0}^{1}|\varphi_{n}(t) (x_{n})|{\rm d}|\nu_{n}|(t)
\rightarrow 0$.\\

(b) If for any $x_{n}\rightarrow 0$ weak it follows that
$\varphi_{n}(t) (x_{n}) \rightarrow 0$ for any $t\in [0,1]${\rm ,}
then $U$ is a Dunford--Pettis operator.\\

(c) $U$ is a compact operator if and only if
$\int_{0}^{1}\|\varphi_{n}(t) \|{\rm d}|\nu_{n}|(t) \rightarrow
0$.\\

\end{enumerate}
\end{theor}

\begin{proof}$\left.\right.$
\begin{enumerate}
\renewcommand\labelenumi{(\roman{enumi})}
\leftskip .2pc
\item Indeed, by hypothesis and the Nikodym boundedness theorem it
follows that $\sup_{n\in \Bbb{N}}\|\lambda_{n}\|([0,1]) <\infty$.
From $\lambda_{n}(E) \rightarrow 0$ $\hbox{weak}^{*}$ for any
$E\in \Sigma$ it follows that for any simple function $f\hbox{\rm
:}\ [0,1] \rightarrow X$ we have $\int_{0}^{1}f(t) {\rm
d}\lambda_{n}(t) \rightarrow 0$. From the well-known inequality
$\|\int_{0}^{1}f(t) {\rm d}\lambda_{n}(t) \|\leq
\|f\|\|\lambda_{n}\|([0,1])$, $f\in B(\Sigma,X)$, we deduce that
$U$ takes its values in $c_{0}$ and that it is bounded linear.

\hskip 1pc The representing measure of $U$ is $G(E)
=(\lambda_{k}(E))_{k\in \Bbb{N}}\hbox{\rm :}\ X\rightarrow c_{0}$
and for any $ x\in X$ we have $\widetilde{G_{x}}([0,1])
=\sup_{E\in \Sigma}\sup_{k\in \Bbb{N}}|\lambda_{k}(E)
(x)|=\sup_{k\in \Bbb{N}}\widetilde{\lambda_{k}x}([0,1])$. Using
that for any $x\in X$ we have $\frac{1}{4}|\lambda_{k}x|([0,1])
\leq \widetilde{\lambda_{k}x}([0,1]) \leq|\lambda_{k}x|([0,1])$.
Applying Theorem~4 we obtain what needs to be proved. (Observe
that in this case, $\|G\|$ is continuous at $\emptyset$ and is
equivalent to the fact that the family $(\lambda_{n})_{n\in
\Bbb{N}}$ is uniformly countably additive (see chapter~I,
Theorem~4, p.~11 of \cite{6}).

\item Let $(\lambda_{n})_{n\in \Bbb{N}}\subseteq
rcabv(\Sigma,X^{*})$ be defined by $\lambda_{n}(E) (x)
=\int_{E}\varphi_{n}(t) (x) \hbox{d}\nu_{n}(t) $ and observe that
$\lambda_{n}(E) = \hbox{Bochner} - \int_{E}\varphi_{n}(t)
\hbox{d}\nu_{n}(t)$.
\end{enumerate}

Then, by hypothesis, we have that $\lambda_{n}(E) \rightarrow 0$
$\hbox{weak}^{*}$ for any $E\in \Sigma$. Also, from
$\|\lambda_{n}(E) \|\leq M|\nu_{n}|(E) $ and the fact that
$(\nu_{n})_{n\in \Bbb{N}}$ is uniformly countably additive it
follows that $(\lambda_{n})_{n\in \Bbb{N}}$ is uniformly countably
additive. Then $U(f) = \big(\int_{0}^{1}f(t)
\hbox{d}\lambda_{n}(t)\big)_{n\in \Bbb{N}}$.

\begin{enumerate}
\renewcommand\labelenumi{(\alph{enumi})}
\leftskip .1pc
\item For any $x\in X$, $k\in \Bbb{N}$ we have that
$\lambda_{k}x\hbox{\rm :}\ \Sigma \rightarrow \Bbb{K}$ is defined
by $(\lambda_{k}x) (E) =\int_{E}\varphi_{k}(t) (x)
\hbox{d}\nu_{k}(t)$ and $|\lambda_{k}x|([0,1])
=\int_{0}^{1}|\varphi_{k}(t) (x)|\hbox{d}|\nu_{k}|(t)$. Now by (i)
if $x_{n}\rightarrow 0$ weak, then $\sup_{k\in
\Bbb{N}}\int_{0}^{1}|\varphi_{k}(t) (x_{n})|\hbox{d}|\nu_{k}|(t)
\rightarrow 0$.

\item Let $(f_{n})_{n\in \Bbb{N}}\subseteq C([0,1],X) $ be such
that $f_{n}\rightarrow 0$ weak. Then $\sup_{n\in \Bbb{N},\, t\in
[0,1]}\|f_{n}(t) \|=L<\infty $ and $f_{n}(t) \rightarrow 0$ weak
for any $t\in [0,1]$.

\hskip 1pc By hypothesis, $\varphi_{n}(t) (f_{n}(t)) \rightarrow
0$ for any $t\in [0,1] $ and $|\varphi_{n}(t) (f_{n}(t))|\leq ML$
for any $t\in [0,1]$. Then, by Bartle's convergence theorem, it
follows that $\int_{0}^{1}|\varphi_{n}(t)
(f_{n}(t))|\hbox{d}|\nu_{n}|(t) \rightarrow 0$ and, by
Corollary~2, $U$ is a Dunford--Pettis operator.

\item We observe that $|\lambda_{n}|([0,1])
=\int_{0}^{1}\|\varphi_{n}(t)\|\hbox{d}|\nu_{n}|(t)$ and the proof
will be finished since (see Exercise~4, p.~114 of \cite{4}) $U$ is
a compact operator if and only if $|\lambda_{n}|([0,1])
\rightarrow 0$.\vspace{-1pc}
\end{enumerate}
\end{proof}

\begin{rema}$\left.\right.$
{\rm
\begin{enumerate}
\renewcommand\labelenumi{(\alph{enumi})}
\leftskip .1pc
\item If $(\varphi_{n})_{n\in \Bbb{N}}\subseteq B(\Sigma,X^{*})$ is
such that $\sup_{n\in \Bbb{N},\,t\in [0,1]}\|\varphi_{n}(t)
\|<\infty${\rm ,} $(\nu_{n})_{n\in \Bbb{N}}\subseteq
rcabv(\Sigma)$ is an uniformly countably additive pointwise
bounded family and in addition{\rm ,} $\varphi_{n}(t) \rightarrow
0$ $\hbox{weak}^{*}$ for any $t\in [0,1]${\rm ,} then it follows
that for any $E\in \Sigma$ and any $x\in X$ we have
$\int_{E}\varphi_{n}(t) (x) \hbox{d}\nu_{n}(t) \rightarrow 0$.

\item Under the hypothesis of theorem~5 (ii){\rm ,} for any $k\in
\Bbb{N}${\rm ,} the operator $T_{k}\hbox{\rm :}\ X\rightarrow
L_{1}(|\nu_{k}|)$ defined by $T_{k}(x) =\varphi_{k}x$ is
Dunford--Pettis and the condition: for any $x_{n}\rightarrow 0$
weak it follows that $\int_{0}^{1}|\varphi_{n}(t)
(x_{n})|\hbox{d}|\nu_{n}|(t) \rightarrow 0$ is equivalent to the
fact that the family $(T_{k})_{k\in \Bbb{N}}$ is an uniformly
Dunford--Pettis family. (For this reason we formulate (a) in
theorem~5(ii)).\vspace{-.5pc}
\end{enumerate}}
\end{rema}

\begin{proof}$\left.\right.$
\begin{enumerate}
\renewcommand\labelenumi{(\alph{enumi})}
\leftskip .1pc
\item Indeed, in our hypotheses, we can apply again Bartle's
convergence theorem, to deduce that for any $E\in \Sigma$ we have
$\int_{E}\varphi_{n}(t) (x) \hbox{d}\nu_{n}(t) \rightarrow
0$.\newpage

\item Let $k\in \Bbb{N}$ be fixed. If $x_{n}\rightarrow 0$ weak,
then $\varphi_{k}(t) (x_{n}) \rightarrow 0$ for any $t\in [0,1]$
and $\sup_{n\in \Bbb{N},\, t\in [0,1]}|\varphi_{k}(t)
(x_{n})|<\infty$. From the Lebesgue dominated convergence theorem,
it follows that $ \int_{0}^{1}|\varphi_{k}(t)
(x_{n})|\hbox{d}|\nu_{k}|(t) \rightarrow 0$. The last part of the
statement follows by Proposition~1.\vspace{-2pc}
\end{enumerate}
\end{proof}

In the following corollary we indicate a way to construct examples
of Dunford--Pettis operators from a Dunford--Pettis one. In view
of Theorem~9 of \cite{2} and Theorem~3.1 of \cite{10}, this result
is, perhaps, natural; see also \cite{9} for other examples in the
scalar case.

\begin{coro}$\left.\right.$
\begin{enumerate}
\renewcommand\labelenumi{\rm (\alph{enumi})}
\leftskip .1pc
\item Let $X$ be a Banach space{\rm ,} $(x_{n}^{*})_{n\in \Bbb{N}}
\subseteq X^{*}$ a bounded sequence{\rm ,} $T\hbox{\rm :}\
X\rightarrow l_{\infty}$ defined by $T(x) =(x_{n}^{*}(x))_{n\in
\Bbb{N}}$ and $(\nu_{n}) _{n\in \Bbb{N}}\subseteq rcabv(\Sigma)$
such that $\nu_{n}(E) \rightarrow 0$ for any $E\in \Sigma$ and
$\liminf_{n\rightarrow \infty}|\nu_{n}|([0,1])
>0${\rm ,} or $x_{n}^{*}\rightarrow 0$ $\hbox{weak}^{*}${\rm ,}
$T\hbox{\rm :}\ X\rightarrow c_{0}$ defined by $T(x) =
(x_{n}^{*}(x))_{n\in \Bbb{N}}$ and $(\nu_{n})_{n\in
\Bbb{N}}\subseteq rcabv(\Sigma)$ uniformly countably additive
pointwise bounded such that $\liminf_{n\rightarrow
\infty}|\nu_{n}|([0,1])>0$.

Let $U\hbox{\rm :}\ C([0,1],X) \rightarrow c_{0}$ be the operator
defined by
\begin{equation*}
\hskip -1.25pc U(f) = \left(\int_{0}^{1}x_{n}^{*}f(t) {\rm
d}\nu_{n}(t) \right)_{n\in \Bbb{N}}.
\end{equation*}
Then\\
%\begin{enumerate}
%\renewcommand\labelenumii{\rm (\roman{enumii})}
\leftskip .2pc (i) $U$ is Dunford--Pettis $\Leftrightarrow$ $T$ is
Dunford--Pettis.\\
(ii) $U$ is compact $\Leftrightarrow$ $T$ is compact.\\

\item Let $a=(a_{n})_{n\in \Bbb{N}}\in l_{\infty}$ and
let $U\hbox{\rm :}\ C([0,1],L_{1}[0,1]) \rightarrow c_{0}$ be the
operator defined by
\begin{equation*}
\hskip -1.25pc U(f)=\left(a_{n}\int_{0}^{1}\left( \int_{0}^{1}f(t)
(s) s^{n}\hbox{d}s \right) r_{n}(t){\rm d}t\right)_{n\in \Bbb{N}}.
\end{equation*}

Then $U$ is Dunford--Pettis{\rm ,} while $U$ is compact
$\Leftrightarrow$ $a\in c_{0}$.
\end{enumerate}
\end{coro}

\begin{proof}$\left.\right.$
\begin{enumerate}
\renewcommand\labelenumi{(\alph{enumi})}
\leftskip .1pc \item Define $\varphi_{n}(t) =x_{n}^{*}$ and
observe that, in our hypotheses, in both cases for any $E\in
\Sigma$ and $ x\in X$ we have $\int_{E}\varphi_{n}(t) (x) {\rm
d}\nu_{n}(t) \rightarrow 0$. Thus, the hypotheses from
Theorem~5(ii) are satisfied and $U(f)=
\big(\int_{0}^{1}\varphi_{n}(t) (f(t)) {\rm
d}\nu_{n}(t)\big)_{n\in \Bbb{ N}}$.\\

%\begin{enumerate}
%\renewcommand\labelenumii{(\roman{enumii})}
\leftskip .2pc (i) Suppose $U$ is a Dunford--Pettis operator. In
both cases, from Theorem~5(ii)(a) if $x_{n}\rightarrow 0$ weak, it
follows that $\int_{0}^{1}|\varphi_{n}(t) (x_{n})|{\rm
d}|\nu_{n}|(t) \rightarrow 0$, or
$|x_{n}^{*}(x_{n})||\nu_{n}|([0,1]) \rightarrow 0$. From here,
since $\liminf_{n\rightarrow \infty}|\nu_{n}|([0,1])>0$, we deduce
$|x_{n}^{*}(x_{n})|\rightarrow 0$ i.e., by Corollary~2, $T$ is
Dunford--Pettis.

\hskip 1pc The converse follows from Corollary~2 and
Theorem~5(ii)(b).

(ii) By Theorem~5(ii)(c), $U$ is compact if and only if
$\|x_{n}^{*}\||\nu_{n}|([0,1]) =\int_{0}^{1}\|\varphi_{n}(t)
\|{\rm d}$ $|\nu_{n}|(t) \rightarrow 0$ or equivalently, by
hypothesis, $\|x_{n}^{*}\|\rightarrow 0$ i.e. $T$ is compact.
%\end{enumerate}\\

We remark that in (i) and (ii) the converses are true without the
hypothesis $\liminf_{n\rightarrow \infty}|\nu_{n}|([0,1])
>0$.

\item Since any positive bounded linear operator from $L_{1}[0,1]$
into $c_{0}$ is Dunford--Pettis (see Corollary~2.3 of \cite{7})
the operator $V\hbox{\rm :}\ L_{1}[0,1] \rightarrow c_{0}$ defined
by $V(f) = \big(\int_{0}^{1}f(s) s^{n}\hbox{d}s\big)_{n\in
\Bbb{N}}$ is Dunford--Pettis, and thus $T\hbox{\rm :}\ L_{1}[0,1]
\rightarrow c_{0}$ defined by $T(f) = \big(a_{n}\int_{0}^{1}f(s)
s^{n}\hbox{d}s\big)_{n\in \Bbb{N}}$ is Dunford--Pettis. By (a)(i),
$U$ is Dunford--Pettis. By (a)(ii), $U$ is compact
$\Leftrightarrow$ $T$ is compact $\Leftrightarrow$
$|a_{n}|\rightarrow 0$.\vspace{-.5pc}
\end{enumerate}

Examples of measures as in Corollary~7(a) can be obtained in the
following ways:
\begin{enumerate}
\renewcommand\labelenumi{\arabic{enumi}.}
\leftskip -.2pc
\item Let $(\alpha_{n})_{n\in \Bbb{N}}\in l_{1}$ be a non-null
element and define $\nu_{n}(E) =\alpha_{1}\int_{E}r_{n}(t)
\hbox{d}t+\dots+\alpha_{n}\int_{E}r_{1}(t) \hbox{d}t$. Then
$\nu_{n}(E) \rightarrow 0$ for any $E\in \Sigma$ and
$\liminf_{n\rightarrow \infty}|\nu_{n}|([0,1]) >0$;

\item If $h\hbox{\rm :}\ [0,1] \rightarrow \Bbb{R}$ is a
differentiable function with $h^{{\prime}}(0) \neq 0$, then
$\nu_{n}(E) = n\int_{E}\big(h\big(\frac{t}{n}\big) -h(0)\big)
\hbox{d}t\rightarrow h^{{\prime}}(0) \int_{E}t\hbox{d}t=\nu (E)$
uniformly for $E\in \Sigma$ and $\lim_{n\rightarrow
\infty}|\nu_{n}|([0,1])
=\frac{|h^{{\prime}}(0)|}{2}>0$.\vspace{-.5pc}
\end{enumerate}

\begin{enumerate}
\renewcommand\labelenumi{\arabic{enumi}.}
\leftskip -.2pc
\item Indeed, in our hypothesis, by a well-known classical result
(see Chapter~IX, Exercise~17 of \cite{11}) it follows that
$\nu_{n}(E) \rightarrow 0$ for any $E\in \Sigma$ and, by
Khinchin's inequality (see p.~10 of \cite{5}) for any $n\in
\Bbb{N}$ we have $|\nu_{n}|([0,1]) \geq \frac{1}{\sqrt{
2}}(|\alpha_{1}|^{2}+\dots +|\alpha_{n}|^{2})^{\frac{1}{2}}$ and
$\liminf_{n\rightarrow \infty}|\nu_{n}|([0,1]) \geq
\frac{1}{\sqrt{2}}\|(\alpha_{n})_{n\in \Bbb{N}}\|_{2}>0$.

\item Let $\varepsilon >0$. Then there is $\delta_{\varepsilon}>0$
such that for any $0\leq t<\delta_{\varepsilon}$ we have
$(h^{{\prime}}(0) -\varepsilon) t\leq h(t) -h(0) \leq
(h^{{\prime}}(0) +\varepsilon) t$. There is also $n_{\varepsilon
}\in \Bbb{N}$ such that $\frac{1}{n_{\varepsilon}}<\delta
_{\varepsilon}$. Take $n\geq n_{\varepsilon}$. Then for any $t\in
[0,1] $ we have $\frac{t}{n}<\delta_{\varepsilon}$ from where
$(h^{{\prime}}(0) -\varepsilon) t\leq n
\big[h\big(\frac{t}{n}\big) -h(0)\big] \leq (h^{{\prime}}(0)
+\varepsilon) t$.

For any $E\in \Sigma$ we obtain $(h^{{\prime}}(0) -\varepsilon)
\int_{E}t\leq n\int_{E} \big[h\big(\frac{t}{n}\big) -h(0)\big]
\hbox{d}t\leq (h^{{\prime}}(0) +\varepsilon) \int_{E}t$, or
$\big|n\int_{E}\big[h\big(\frac{t}{n}\big)-h(0)\big]
\hbox{d}t-h^{{\prime}}(0) \int_{E}t\hbox{d}t\big|\leq \varepsilon
\int_{E}t\hbox{d}t\leq \frac{\varepsilon}{2}$. Also
\begin{align*}
\hskip -1.25pc |\nu_{n}|([0,1]) = n\int_{0}^{1} \left|h \left(\!
\frac{t}{n} \!\right) -h(0) \right|\hbox{d}t = n^{2}
\int_{0}^{\frac{1}{n}}|h(t)\! -h(0)|\hbox{d}t\rightarrow
\frac{|h^{{\prime}}(0)|}{2}.
\end{align*}

The next example is different from what was used in Theorem~9 of
\cite{2} and Theorem~3.1 of \cite{10} and in the scalar case
appear in Example~11 of \cite{9}.\vspace{-1pc}
\end{enumerate}
\end{proof}

\begin{coro}$\left.\right.$\vspace{.5pc}

\noindent Let $X$ be a Banach space{\rm ,}
$\sum_{n=1}^{\infty}x_{n}^{*}$ a weakly Cauchy series in $X^{*}$
and let $U\hbox{\rm :}\ C([0,1],X) \rightarrow c_{0}$ be the
operator defined by
\begin{equation*}
U(f) = \left( \int_{0}^{1}(x_{n}^{*}(f(t)) r_{1}(t) +\dots
+x_{1}^{*}(f(t)) r_{n}(t)) \hbox{\rm d}t\right)_{n\in \Bbb{N}}.
\end{equation*}

Then $U$ is a Dunford--Pettis operator and $U$ is a compact
operator if and only if $U=0$.
\end{coro}

\begin{proof}
Let $T\hbox{\rm :}\ X\rightarrow l_{1}$ be the operator defined by
$T(x) =(x_{n}^{*}(x))_{n\in \Bbb{N}}$ and $V\hbox{\rm :}\
C([0,1],X) \rightarrow C([0,1],l_{1}) $ the operator defined by
$V(f) =T\circ f$. Define also $S\hbox{\rm :}\ C([0,1],l_{1})
\rightarrow c_{0}$ by
\begin{equation*}
S(f) = \left(\int_{0}^{1}(\langle f(t),e_{n}\rangle r_{1}(t) +
\dots + \langle f(t),e_{1}\rangle r_{n}(t)) \hbox{d}t
\right)_{n\in \Bbb{N}}
\end{equation*}
and observe that $U=SV$.\newpage

($S$ takes its values in $c_{0}$ (Chapter~IX, Exercise~17 of
\cite{11}).) Because $l_{1}$ has the Schur property, from
Theorem~3.1 of \cite{10}, $S$ is a Dunford--Pettis operator and
hence $U$ is also a Dunford--Pettis operator.

Define $\varphi_{n}(t) =x_{n}^{*}r_{1}(t) +x_{n-1}^{*}r_{2}(t) +
\dots + x_{1}^{*}r_{n}(t) $ and observe that (Chapter~IX,
Exercise~17 of \cite{11}) $\int_{E}\varphi_{n}(t)
\hbox{d}t\rightarrow 0$ $\hbox{weak}^{*}$ for any $E\in \Sigma$,
$\|\varphi _{n}(t) \|\leq w_{1}(x_{n}^{*}| n\in \Bbb{N})$ for any
$t\in [0,1]$ and $U(f) = \big(\int_{0}^{1}\varphi_{n}(t) (f(t))
\hbox{d}t\big)_{n\in \Bbb{N}}$.

By Theorem~5(ii)(c) $U$ is compact if and only if
$\int_{0}^{1}\|\varphi_{n}(t) \|\hbox{d}t\rightarrow 0$. However,
for any $ n\in \Bbb{N}$, again by Khinchin's inequality we have
\begin{equation*}
\frac{1}{\sqrt{2}}\sup\limits_{\|x\|\leq 1}(|x_{1}^{*}(x)|^{2} +
\dots + |x_{n}^{*}(x)|^{2})^{\frac{1}{2}}\leq \int_{0}^{1}\|
\varphi_{n}(t) \|\hbox{d}t.
\end{equation*}

If $U$ is compact, then $\sup_{\|x\|\leq 1}
\big(\sum_{n=1}^{\infty} |x_{n}^{*}(x)|^{2}\big)^{\frac{1}{2}}=0$,
which implies $x_{n}^{*}=0$ for any $n\in \Bbb{N}$ i.e. $U=0$.
\end{proof}

We now state a remark which is certainly well-known, but,
unfortunately, we do not know a reference.

\begin{rema}$\left.\right.$

{\rm
\begin{enumerate}
\renewcommand\labelenumi{\rm (\alph{enumi})}
\leftskip .1pc
\item The space $w_{1}(L_{1}[0,1]^{*})$ is isometrically
isomorph with $L_{\infty}([0,1],l_{1})$, more precisely, if
$g_{n}\in L_{\infty}[0,1] =L_{1}[0,1] ^{*}$, then
$\sum_{n=1}^{\infty}g_{n}$ is a weakly Cauchy series in
$L_{\infty}[0,1]$ $\Leftrightarrow$ the function $g=(g_{n})_{n\in
\Bbb{N}}\in L_{\infty}([0,1],l_{1})$.

\item A~weakly Cauchy series $\sum_{n=1}^{\infty}g_{n}$
in $L_{\infty}[0,1]$ is unconditionally norm convergent
$\Leftrightarrow$ the function $g=(g_{n})_{n\in \Bbb{N}}\in
L_{\infty}([0,1],l_{1})$ has an essentially relatively compact
range in $l_{1}$.\vspace{-.5pc}
\end{enumerate}}
\end{rema}

\begin{proof}$\left.\right.$

\begin{enumerate}
\renewcommand\labelenumi{\rm (\alph{enumi})}
\leftskip .1pc
\item Indeed, $\sum_{n=1}^{\infty}g_{n}$ is a
weakly Cauchy series in $L_{\infty}[0,1]$ $\Leftrightarrow$ the
operator $T\hbox{\rm :}\ L_{1}[0,1] \rightarrow l_{1}$ defined by
$T(f) = \big(\int_{0}^{1}g_{n}(t) f(t) \hbox{d}t\big)_{n\in
\Bbb{N}}$ is bounded linear. Since $l_{1}$ has the Radon--Nikodym
property (see Theorem, p.~63 of \cite{6}), this is equivalent to
the fact that $T$ is representable i.e. there is $h=(h_{n})_{n\in
\Bbb{N}}\in L_{\infty}([0,1],l_{1})$ such that $T(f)
=\hbox{Bochner}-\int_{0}^{1}f(t) h(t)\hbox{d}t$. Then for any
$n\in \Bbb{N}$ we have that $\int_{E}g_{n}(t)
\hbox{d}t=\int_{E}h_{n}(t)\hbox{d}t$ for any $E\in \Sigma$. Thus
$g_{n}=h_{n}$ $\mu$-a.e. and the statement follows.

\item By Theorem~1.9, p.~9 of \cite{5}, the unconditionality norm
convergence of series $\sum_{n=1}^{\infty}g_{n}$ is equivalent to
the fact that the operator $T\hbox{\rm :}\ L_{1}[0,1] \rightarrow
l_{1}$ as in (a) is compact. By the representation of compact
operators on $L_{1}(\mu)$ (see p.~68 of \cite{6}), this is
equivalent to the fact that $g$ has an essentially relatively
compact range in $l_{1}$.\vspace{-2pc}
\end{enumerate}
\end{proof}

In the sequel we analyze the same kind of operators as in
Corollary~7, but with values in $l_{p}$, where $1\leq p<\infty$.
We begin with a lemma which is, probably, a well-known result but,
we do not know a reference.

\begin{lem} Let $1\leq p<\infty${\rm ,} $(\alpha_{n})_{n\in
\Bbb{N}}\in l_{p}$ and let $G\hbox{\rm :}\ \Sigma \rightarrow
l_{p}$ be defined by
\begin{equation*}
G(E) = \left(\alpha_{n}\int_{E}r_{n}(t) \hbox{\rm d}t
\right)_{n\in \Bbb{N}}.
\end{equation*}

\begin{enumerate}
\renewcommand\labelenumi{\rm (\roman{enumi})}
\leftskip .2pc
\item If $1\leq p<2${\rm ,} then
$\frac{1}{\sqrt{2}}\|(\alpha_{n})_{n\in \Bbb{N}}\|_{r}\leq
\|G\|([0,1]) \leq \|(\alpha_{n})_{n\in \Bbb{N}}\|_{r}${\rm ,}
where $\frac{1}{p}=\frac{1}{2}+\frac{1}{r}$.

\item If $2\leq p<\infty${\rm ,} then
$\frac{1}{\sqrt{2}}\sup_{n\in \Bbb{N}}|\alpha_{n}|\leq
\|G\|([0,1]) \leq \sup_{n\in \Bbb{N}}|\alpha_{n}|$.\vspace{-.5pc}
\end{enumerate}
\end{lem}

\begin{proof}
Let $h\hbox{\rm :}\ [0,1] \rightarrow l_{p}$ be defined by $h(t)
=(\alpha_{n}r_{n}(t))_{n\in \Bbb{ N}}$ and observe that $G(E)
=\hbox{Bochner}-\int_{E}h(t) \hbox{d}t$. Then
\begin{equation*}
\|G\|([0,1]) =\sup\limits_{\|y^{*}\|\leq 1}|G_{y^{*}}|([0,1])
=\sup\limits_{\|y^{*}\|\leq 1}\int_{0}^{1}|y^{*}h(t)|\hbox{d}t,
\end{equation*}
because $G_{y^{*}}(E) = \int_{E}y^{*}h(t) \hbox{d}t$ and
$|G_{y^{*}}|([0,1]) =\int_{0}^{1}|y^{*}h(t)|\hbox{d}t$.

However, for any $y^{*}=(\xi_{n})_{n\in \Bbb{N}}\in
l_{p}^{*}=l_{q}$ ($q$ is the conjugate of $p$), we have $y^{*}h(t)
= \sum_{n=1}^{\infty}\xi_{n}\alpha_{n}r_{n}(t)$ and, by Khinchin's
inequality
\begin{equation*}
\frac{1}{\sqrt{2}} \left( \sum\limits_{n=1}^{\infty}|\xi_{n}
\alpha_{n}|^{2} \right)^{\frac{1}{2}}\leq \int_{0}^{1} |y^{*}h(t)|
\hbox{d}t\leq \left(\sum\limits_{n=1}^{\infty}
|\xi_{n}\alpha_{n}|^{2} \right)^{\frac{1}{2}},
\end{equation*}
i.e. $\frac{1}{\sqrt{2}}\|M\|\leq \|G\|([0,1]) \leq \|M\|$, where
$ M\hbox{\rm :}\ l_{q}\rightarrow l_{2}$ is the multiplication
operator $M((\xi_{n})_{n\in \Bbb{N}}) =(\alpha_{n}\xi_{n})_{n\in
\Bbb{N}}$. Now, as is well-known

\begin{enumerate}
\renewcommand\labelenumi{\rm (\roman{enumi})}
\leftskip .2pc
\item
if $2<q$ i.e. $1\leq p<2$, then $\|M\|=\|(\alpha_{n})_{n\in
\Bbb{N}}\|_{r}= \big(\sum_{n=1}^{\infty} |\alpha_{n}|^{r}
\big)^{\frac{1}{r}}$, where $\frac{1}{2}=\frac{1}{q}+\frac{1}{r}$
i.e. $\frac{1}{p}=\frac{1}{2}+\frac{1}{r}$.

\item if $q\leq 2$ i.e. $2\leq p<\infty$, then
$\|M\|=\sup_{n\in \Bbb{N}}|\alpha_{n}|$.\vspace{-1.8pc}
\end{enumerate}
\end{proof}

In case of operators on $C([0,1],X) $\ with values in $l_{p}$,
Khinchin's inequality gives a distinction for $1<p<2$ and $2\leq
p<\infty$.

\begin{propo}$\left.\right.$\vspace{.5pc}

\noindent Let $1<p<\infty${\rm ,} $p^{*}$ be the conjugate of
$p${\rm ,} $X$ a Banach space{\rm ,} $(x_{n}^{*})_{n\in
\Bbb{N}}\in w_{p}(X^{*})$ and $T\hbox{\rm :}\ X\rightarrow l_{p}$
defined by $T(x) =(x_{n}^{*}(x))_{n\in \Bbb{N}}$.

Let $U\hbox{\rm :}\ C([0,1],X) \rightarrow l_{p}$ be the operator
defined by
\begin{equation*}
U(f) = \left( \int_{0}^{1}x_{n}^{*}f(t) r_{n}(t){\rm d}t
\right)_{n\in \Bbb{N}}.
\end{equation*}

\begin{enumerate}
\renewcommand\labelenumi{\rm (\alph{enumi})}
\leftskip .1pc
\item If $T$ is a Dunford--Pettis {\rm (}resp. compact{\rm
)} operator{\rm ,} then $U$ is a Dunford--Pettis {\rm (}resp.
compact{\rm )} operator.

\item If $X$ is reflexive and $T$ is Dunford--Pettis{\rm ,}
or $X=c_{0}${\rm ,} then $T$ is compact and thus $U$ is compact.

\item If $U$ is a compact operator{\rm ,} then
$x_{n}^{*}\rightarrow 0$ in norm.

\item Suppose $1<p<2$. Then

%\begin{enumerate}
%\renewcommand\labelenumii{\rm (\roman{enumii})}
\leftskip .2pc (i) $U$ is a Dunford--Pettis operator
$\Leftrightarrow$ $T$ is a Dunford--Pettis operator.

(ii) $U$ is a compact operator $\Leftrightarrow$ $T$ is a compact
operator.
%\end{enumerate}

\item Suppose $2\leq p<\infty$.

%\begin{enumerate}
%\renewcommand\labelenumii{\rm (\roman{enumii})}
\leftskip .2pc (i) If $U$ is a Dunford--Pettis operator{\rm ,}
then $T\hbox{\rm :}\ X\rightarrow c_{0}$ is a Dunford--Pettis
operator.

(ii) If $T\hbox{\rm :}\ X\rightarrow c_{0}$ is a compact operator
i.e. $x_{n}^{*}\rightarrow 0$ in norm and $X^{*}$ has type $a${\rm
,} where $1<p^{*}\leq a\leq 2${\rm ,} then $U$ is a compact
operator.\vspace{-1pc}
%\end{enumerate}
\end{enumerate}
\end{propo}

\begin{proof} The fact that $U$ is well-defined, bounded linear is
clear. Also, the representing measure of $U$ is $G(E) (x)
=\big(x_{n}^{*}(x) \int_{E}r_{n}(t) \hbox{d}t\big)_{n\in
\Bbb{N}}$.

\noindent (a) We observe that for any $f\in C([0,1],X)$ we have
\begin{align*}
\|U(f) \|^{p} &=
\sum\limits_{n=1}^{\infty}\left|\int_{0}^{1}x_{n}^{*}f(t)
r_{n}(t)\hbox{d}t\right|^{p}\leq
\int_{0}^{1}\sum\limits_{n=1}^{\infty}
|x_{n}^{*}f(t)|^{p}\hbox{d}t\\[.4pc]
&= \int_{0}^{1}\|T(f(t)) \|^{p}\hbox{d}t
\end{align*}
i.e.
\begin{equation*}
\hskip -4pc (*)\hskip 2.9pc \|U(f) \|\leq \left(
\int_{0}^{1}\|T(f(t))\|^{p}\hbox{d}t \right)^{1/p}.
\end{equation*}

From this inequality it is easy to prove that if $T$ is a
Dunford--Pettis operator, then $U$ is a Dunford--Pettis operator.

Indeed, let $(f_{n})_{n\in \Bbb{N}}\subseteq C([0,1],X) $ be such
that $f_{n}\rightarrow 0$ weak. Then $\sup_{n\in \Bbb{N},\, t\in
[0,1]}$ $\|f_{n}(t) \|=L<\infty $ and $f_{n}(t) \rightarrow 0$
weak for any $t\in [0,1]$. Since $T$ is a Dunford--Pettis
operator, it follows that $\|T(f_{k}(t)) \|^{p}\rightarrow 0$ for
any $t\in [0,1]$ and $\|T(f_{k}(t)) \|\leq \|T\|\|f_{k}(t) \|\leq
L\|T\|$.

From the Lebesgue dominated convergence theorem we get
$\int_{0}^{1}\|T(f_{k}(t)) \|^{p}\hbox{d}t\rightarrow 0$ and by
$(*)$, it follows that $\|U(f_{k}) \|\rightarrow 0$.

Suppose $T\hbox{\rm :}\ X\rightarrow l_{p}$ is compact and let
$\varepsilon >0$. From Exercise~6, p.~6 of \cite{4}, there is
$n_{\varepsilon}\in \Bbb{N}$ such that $\sup_{\|x\|\leq
1}\sum_{k=n_{\varepsilon}}^{\infty}|x_{k}^{*}(x)|^{p}<\varepsilon$.

Take $f\in C([0,1],X) $ with $\|f\|\leq 1$. Since for any $n\in
\Bbb{N}$ we have $\sum_{k=n}^{\infty}
\big|\int_{0}^{1}x_{n}^{*}f(t)$ $r_{n}(t)\hbox{d}t\big|^{p}\leq
\int_{0}^{1}$ $\sum_{k=n}^{\infty}|x_{k}^{*}(f(t))|^{p}\hbox{d}t$,
we deduce $\sum_{k=n_{\varepsilon}}^{\infty}
\big|\int_{0}^{1}x_{n}^{*}f(t) r_{n}(t)\hbox{d}t\big|^{p}\leq
\varepsilon$, hence, again by Exercise~6, p.~6 of \cite{4}, $U$ is
compact.

\noindent (b) First, it follows from the well-known fact that a
Dunford--Pettis operator whose domain is reflexive is compact and
second, by the same reasoning and the fact that $l_{1}$ has the
Schur property, the dual is compact, hence compact by Schauder's
theorem.

\noindent (c) We have that $U^{*}\hbox{\rm :}\
l_{p^{*}}\rightarrow rcabv(\Sigma,X^{*}) =C([0,1],X)^{*}$ is
defined by $U^{*}(\xi) =G_{\xi}$, where $G_{\xi}(E) (x) =
\sum_{n=1}^{\infty} \xi_{n}x_{n}^{*}(x) \int_{E}r_{n}(t)
\hbox{d}t$.

Since $\xi =(\xi_{n})_{n\in \Bbb{N}}\in l_{p^{*}}$ and
$(x_{n}^{*})_{n\in \Bbb{N}}\in w_{p}(X^{*}) $ it follows that the
series $\sum_{n=1}^{\infty}\xi_{n}x_{n}^{*}r_{n}(t) $ is
unconditionally norm convergent for any $t\in [0,1] $ and let
$h_{\xi}\hbox{\rm :}\ [0,1] \rightarrow X^{*}$ be defined by
$h_{\xi}(t) =\sum_{n=1}^{\infty}\xi_{n}x_{n}^{*}r_{n}(t)$. Then
the function $h_{\xi}$ is Bochner integrable (by an easy
application of the Lebesgue dominated convergence theorem for the
Bochner integral), $G_{\xi}(E) = \hbox{Bochner}-\int_{E}h_{\xi}(t)
\hbox{d}t$ and thus $\|U^{*}(\xi) \|=|G_{\xi}|([0,1])
=\int_{0}^{1}\|h_{\xi }(t) \|\hbox{d}t$.

If $U$ is a compact operator, then by the Schauder's theorem
$U^{*}$ is compact, in particular $U^{*}(e_{n}) \rightarrow 0$ and
the statement follows since $\|U^{*}(e_{n})
\|=\int_{0}^{1}\|h_{e_{n}}(t) \|\hbox{d}t= \int_{0}^{1}
\|x_{n}^{*}r_{n}(t) \|\hbox{d}t=\|x_{n}^{*}\|$.

\noindent (d) Define $r$ such that
$\frac{1}{p}=\frac{1}{2}+\frac{1}{r}$ and $T\hbox{\rm :}\
X\rightarrow l_{r}$ by $T(x) =(x_{n}^{*}(x))_{n\in \Bbb{N}}$.

We prove the following equivalences:

\begin{enumerate}
\renewcommand\labelenumi{\rm (\roman{enumi})}
\leftskip .2pc
\item $U$ is a Dunford--Pettis operator $\Leftrightarrow$ $T\hbox{\rm :}\
X\rightarrow l_{r}$ is a Dunford--Pettis operator
$\Leftrightarrow$ $T\hbox{\rm :}\ X\rightarrow l_{p}$ is a
Dunford--Pettis operator.

\hskip 1pc We observe that for any $x\in X$ the measure
$G_{x}\hbox{\rm :}\ \Sigma \rightarrow l_{p}$ is defined by
$G_{x}(E) = \big(x_{n}^{*}(x) \int_{E}r_{n}(t)
\hbox{d}t\big)_{n\in \Bbb{N}}$ and by lemma~10(i)
\begin{equation*}
\hskip -1.25pc \frac{1}{\sqrt{2}}\|T(x) \|_{r}\leq
\|G_{x}\|([0,1]) \leq \|T(x) \|_{r}.
\end{equation*}

If $U$ is a Dunford--Pettis operator, then by Theorem~4, for any
$x_{n}\rightarrow 0$ weak it follows that $\|G_{x_{n}}\|([0,1])
\rightarrow 0$ and by the above $\|T(x_{n}) \|_{r}\rightarrow 0$
i.e. $T\hbox{\rm :}\ X\rightarrow l_{r}$ is Dunford--Pettis.

\hskip 1pc We now prove that for any $x\in X$ we have the
inequalities
\begin{equation*}
\hskip -3.8pc (**) \hskip 1.1pc \|T(x) \|_{r}\leq \|T(x)
\|_{p}\leq w_{2}(x_{n}^{*}| n\in \Bbb{N}) \|x\|\|T(x) \|_{r}
\end{equation*}
and then from these inequalities it is easy to prove the second
equivalence.

\hskip 1pc Indeed, first inequality follows from the inclusion
$l_{p}\subset l_{r}$ and for the second we use the H\"{o}lder
inequality $\|T(x) \|_{p}\leq \|T(x) \|_{2}\|T(x) \|_{r}$.

\hskip 1pc In (a) we prove that if $T\hbox{\rm :}\ X\rightarrow
l_{p}$ is a Dunford--Pettis operator, then $U$ is a
Dunford--Pettis operator.

\item $U$ is a compact operator $\Leftrightarrow$ $T\hbox{\rm :}\
X\rightarrow l_{r}$ is a compact operator $\Leftrightarrow$
$T\hbox{\rm :}\ X\rightarrow l_{p}$ is a compact operator.

\hskip 1pc Suppose $U$ is compact and $\varepsilon >0$. Then (see
Theorem~6 of \cite{1} and Exercise~6, p.~6 of \cite{4}) there is
$n_{\varepsilon}\in \Bbb{N}$ such that for any $f\in B(\Sigma,X)$
with $\|f\|\leq 1$ we have
\begin{equation*}
\hskip -1.25pc \left(\sum\limits_{k=n_{\varepsilon}}^{\infty}
\left| \int_{0}^{1}x_{k}^{*}f(t) r_{k}(t)\hbox{d}t \right|^{p}
\right)^{1/p}<\varepsilon.
\end{equation*}

By H\"{o}lder's inequality we obtain that for any $f\in
B(\Sigma,X)$ with $\|f\|\leq 1$ and any $\xi =(\xi_{n})_{n\in
\Bbb{N}}\in l_{p^{*}}$ we have
\begin{equation*}
\hskip -1.25pc \left|\sum_{k=n_{\varepsilon}}^{\infty}
\xi_{k}\int_{0}^{1} x_{k}^{*}f(t) r_{k}(t)\hbox{d}t \right|\leq
\varepsilon \left( \sum\limits_{k=n_{\varepsilon}}^{\infty}
|\xi_{k}|^{p^{*}} \right)^{1/p^{*}}.
\end{equation*}

In particular, for any $E\in \Sigma$ and $\|x\|\leq 1$ and any
$\xi \in l_{p^{*}}$ we have
\begin{equation*}
\hskip -1.25pc \left| \int_{E}
\left(\sum\limits_{k=n_{\varepsilon}}^{\infty} \xi_{k}x_{k}^{*}(x)
r_{k}(t) \right) \hbox{d}t \right|\leq \varepsilon \left(
\sum\limits_{k=n_{\varepsilon}}^{\infty} |\xi_{k}|^{p^{*}}
\right)^{1/p^{*}}.
\end{equation*}
Then for any $\|x\|\leq 1$ and any $\xi \in l_{p^{*}}$ we deduce
\begin{align*}
\hskip -1.25pc \frac{1}{4}\int_{0}^{1} \left|
\sum\limits_{k=n_{\varepsilon}}^{\infty} \xi_{k}x_{k}^{*}(x)
r_{k}(t) \right| \hbox{d}t &\leq \sup\limits_{E\in \Sigma} \left|
\int_{E} \left( \sum\limits_{k=n_{\varepsilon}}^{\infty}
x_{k}^{*}(x) r_{k}(t) \right) \hbox{d}t\right|\\[.4pc]
\hskip -1.25pc &\leq \varepsilon \left(
\sum\limits_{k=n_{\varepsilon}}^{\infty} |\xi_{k}|^{p^{*}}
\right)^{1/p^{*}}
\end{align*}
and, by Khinchin's inequality we get
\begin{equation*}
\hskip -1.25pc \frac{1}{4\sqrt{2}} \left(
\sum\limits_{k=n_{\varepsilon}}^{\infty} |\xi_{k}x_{k}^{*}(x)|^{2}
\right)^{1/2} \leq \varepsilon \left(
\sum\limits_{k=n_{\varepsilon}}^{\infty} |\xi_{k}|^{p^{*}}
\right)^{1/p^{*}}.
\end{equation*}

Taking $\xi_{k}\in \Bbb{K}$ such that
$|\xi_{k}|=|x_{k}^{*}(x)|^{(r-2)/2}$ and using
$\frac{1}{p^{*}}=\frac{1}{2}-\frac{1}{r}$ we obtain
$\frac{1}{4\sqrt{2}}\big(\sum_{k=n_{\varepsilon}}^{\infty}
|x_{k}^{*}(x)|^{r} \big)^{1/r}\leq \varepsilon $ i.e., by
Exercise~6, p.~6 of \cite{4}, $T\hbox{\rm :}\ X\rightarrow l_{r}$
is compact.

Now, from the inequality $(**)$ in (i) it can be proved, using
Exercise~6, p.~6 of \cite{4}, that if $T\hbox{\rm :}\ X\rightarrow
l_{r}$ is a compact operator, then $T\hbox{\rm :}\ X\rightarrow
l_{p}$ is a compact operator and by (a), $U$ is
compact.\vspace{-1pc}
\end{enumerate}

\noindent (e)\vspace{-.5pc}

\begin{enumerate}
\renewcommand\labelenumi{\rm (\roman{enumi})}
\leftskip .2pc
\item Is the same as in (d)(i) and use Lemma~10(ii). We omit the
proof. Another way to prove this fact is to use Theorem~5(ii)(a).

\item We prove that in our hypothesis it follows that $U^{*}$ is
compact, hence by Schauder's theorem $U$ will be compact. With the
same notations as in (c), $U^{*}\hbox{\rm :}\ l_{p^{*}}\rightarrow
rcabv(\Sigma,X^{*}) $ is defined by $U^{*}(\xi) =G_{\xi}$ and
\begin{equation*}
\hskip -1.25pc \|U^{*}(\xi) \|=|G_{\xi}|([0,1])
=\int_{0}^{1}\|h_{\xi}(t) \|\hbox{d}t.
\end{equation*}

Since $X^{*}$ has type $a$, then by definition of type $a$ (see
chapter~11, p.~217 of \cite{5}) for any $\xi =(\xi_{n})_{n\in
\Bbb{N}}\in l_{p^{*}}$, we have
\begin{equation*}
\hskip -1.25pc \|U^{*}(\xi)\| = \int_{0}^{1}\|h_{\xi}(t)
\|\hbox{d}t\leq T_{a}(X^{*}) \left( \sum\limits_{n=1}^{\infty}
|\xi_{n}|^{a} \|x_{n}^{*}\|^{a} \right)^{1/a}.
\end{equation*}

From $p^{*}\leq a$ we obtain that
\begin{equation*}
\hskip -1.25pc \|U^{*}(\xi) \|\leq T_{a}(X^{*}) \left(
\sum\limits_{n=1}^{\infty} |\xi_{n}|^{p^{*}}\|x_{n}^{*}\|^{p^{*}}
\right)^{1/p^{*}}.
\end{equation*}

If we define the finite rank operators $V_{n}\hbox{\rm :}\
l_{p^{*}}\rightarrow rcabv(\Sigma,X^{*})$ by $V_{n}(\xi)
=U^{*}(\xi_{1},\dots,\xi_{n},0,\dots)$, then by the above
inequality, for any $\xi =(\xi_{n})_{n\in \Bbb{N}}\in l_{p^{*}}$
with $\|\xi \|\leq 1$ we obtain
\begin{align*}
\hskip -1.25pc \|U^{*}(\xi) -V_{n}(\xi)\| &\leq T_{a}(X^{*})
\left( \sum\limits_{k=n+1}^{\infty} |\xi_{k}|^{p^{*}}
\|x_{k}^{*}\|^{p^{*}} \right)^{1/p^{*}}\\[.4pc]
\hskip -1.25pc &\leq T_{a}(X^{*}) \sup\limits_{k\geq
n+1}\|x_{k}^{*}\|
\end{align*}
and hence
\begin{equation*}
\hskip -1.25pc \|U^{*}-V_{n}\|\leq T_{a}(X^{*}) \sup\limits_{k\geq
n+1}\|x_{k}^{*}\|\rightarrow 0.
\end{equation*}
Thus $U^{*}$ is compact.
\end{enumerate}
\end{proof}
With the help of Proposition~11 we can give\vspace{0.2cm} some
concrete examples.

\begin{exam}$\left.\right.$
{\rm
\begin{enumerate}
\renewcommand\labelenumi{\rm (\alph{enumi})}
\leftskip .15pc
\item The operator $U\hbox{\rm :}\
C([0,1],L_{2}[0,1]) \rightarrow l_{2}$ defined by
\begin{equation*}
\hskip -1.25pc U(f)= \left( \int_{0}^{1} \left( \int_{0}^{1}f(t)
(s) s^{n}\hbox{d}s \right) r_{n}(t)\hbox{d}t \right)_{n\in
\Bbb{N}}
\end{equation*}
is compact.

\item Let $1<p<\infty$ and $U\hbox{\rm :}\ C([0,1],l_{p})
\rightarrow l_{p}$ be the operator defined by
\begin{equation*}
\hskip -1.25pc U(f) = \left( \frac{\int_{0}^{1}(\langle
f(t),e_{1}\rangle + \dots + \langle f(t),e_{n}\rangle) r_{n}(t)
\hbox{d}t}{n}\right)_{n\in \Bbb{N}}.
\end{equation*}

%\begin{enumerate}
%\renewcommand\labelenumii{\rm (\roman{enumii})}
\leftskip .4pc (i) If $1<p<2$, then $U$ is not Dunford--Pettis.

(ii) If $2\leq p<\infty${\rm ,} then $U$ is compact.
%\end{enumerate}

\item Let $\alpha =(\alpha_{n})_{n\in \Bbb{N}}\in l_{\infty}$,
$1<p<\infty$ and $U\hbox{\rm :}\ C([0,1],l_{1}) \rightarrow l_{p}$
be the operator defined by
\begin{equation*}
\hskip -1.25pc U(f) = \left( \alpha_{n}\int_{0}^{1}\langle f(t),
e_{n}\rangle r_{n}(t)\hbox{d}t \right)_{n\in \Bbb{N}}.
\end{equation*}
Then $U$ is a Dunford--Pettis operator, and $U$ is compact
$\Leftrightarrow$ $\alpha \in c_{0}$.

\item Let $1<p<2$. Define $\frac{1}{p}=\frac{1}{2}+\frac{1}{r}$
and take $\alpha =(\alpha_{n})_{n\in \Bbb{N}}\in l_{r}$. Then the
operator $U\hbox{\rm :}\ C([0,1],C[0,1]) \rightarrow l_{p}$
defined by
\begin{equation*}
\hskip -1.25pc U(f)= \left(\alpha_{n}\int_{0}^{1}
\left(\int_{0}^{1}f(t) (s) r_{n}(s) \hbox{d}s\right)
r_{n}(t)\hbox{d}t\right)_{n\in \Bbb{N}}
\end{equation*}
is compact.

\item Let $2\leq p<\infty$, $\alpha =(\alpha_{n})_{n\in
\Bbb{N}}\in l_{\infty}$ and let $U\hbox{\rm :}\ C([0,1],C[0,1])
\rightarrow l_{p}$ be the operator defined by
\begin{equation*}
\hskip -1.25pc U(f) = \left( \alpha_{n}\int_{0}^{1}
\left(\int_{0}^{1}f(t) (s) r_{n}(s) \hbox{d}s \right)
r_{n}(t)\hbox{d}t \right)_{n\in \Bbb{N}}.
\end{equation*}
Then $U$ is a Dunford--Pettis operator, and $U$ is compact
$\Leftrightarrow$ $\alpha \in c_{0}$.

\item Let $2\leq p<\infty$, $1<p^{*}\leq 2$ be the
conjugate of $p$, $\alpha =(\alpha_{n})_{n\in \Bbb{N}}\in
l_{\infty}$ and $U\hbox{\rm :}\ C([0,1],L_{p^{*}}[0,1])
\rightarrow l_{p}$ the operator defined by
\begin{equation*}
\hskip -1.25pc U(f)= \left( \alpha_{n}\int_{0}^{1} \left(
\int_{0}^{1}f(t) (s) r_{n}(s) \hbox{d}s \right) r_{n}(t)\hbox{d}t
\right)_{n\in \Bbb{N}}.
\end{equation*}
Then $U$ is a Dunford--Pettis operator $\Leftrightarrow$ $U$ is
compact $\Leftrightarrow$ $\alpha \in c_{0}$.

\item Let $1<p<2$, $g=(g_{n})_{n\in \Bbb{N}}\in L_{\infty}([0,1],
l_{p})$ and let $U\hbox{\rm :}\ C([0,1],L_{1}[0,1]) \rightarrow
l_{p}$ be the operator defined by
\begin{equation*}
\hskip -1.25pc U(f) = \left( \int_{0}^{1} \left( \int_{0}^{1}f(t)
(s) g_{n}(s) \hbox{d}s \right) r_{n}(t)\hbox{d}t \right)_{n\in
\Bbb{N}}.
\end{equation*}
Then $U$ is a Dunford--Pettis operator, and $U$ is compact
$\Leftrightarrow$ $g$ has an essentially relatively compact
range.\vspace{-.5pc}
\end{enumerate}}
\end{exam}

\begin{proof}$\left.\right.$

\begin{enumerate}
\renewcommand\labelenumi{\rm (\alph{enumi})}
\leftskip .15pc
\item Take $x_{n}^{*}\in L_{2}[0,1]$ defined by $x_{n}^{*}(f) =
\int_{0}^{1}f(s) s^{n}\hbox{d}s$. By Hilbert's theorem, the
operator $H\hbox{\rm :}\ L_{2}[0,1] \rightarrow l_{2}$ defined by
$H(f) = \big(\int_{0}^{1}f(s) s^{n}\hbox{d}t\big)_{n\in \Bbb{N}}$
is bounded linear. Now $H\hbox{\rm :}\ L_{2}[0,1] \rightarrow
c_{0}$ is compact, because $\|x_{n}^{*}\|=\frac{1}{\sqrt{2n}}
\rightarrow 0$, $L_{2}[0,1]$ is a Hilbert space, hence of type
$2$, and by Proposition~11(e)(ii) we obtain the statement.

\item Let $x_{n}^{*}\in l_{p}^{*}$ be defined by
$x_{n}^{*}((x_{n})_{n\in \Bbb{N}}) = \frac{x_{1}+ \dots +
x_{n}}{n}$. Then the Hardy operator $H\hbox{\rm :}\
l_{p}\rightarrow l_{p}$ defined by $H(x) =(x_{n}^{*}(x))_{n\in
\Bbb{N}}$ is bounded linear.

%\begin{enumerate}
%\renewcommand\labelenumii{\rm (\roman{enumii})}
\leftskip .4pc (i) It follows from Proposition~11(d)(i) and (b)
and the well-known fact that the Hardy operator is not compact.

(ii) In this case $l_{p}^{*}=l_{p^{*}}$ has type $p^{*}$, where
$1<p^{*}\leq 2$, $\|x_{n}^{*}\| =
\frac{1}{n^{\frac{1}{p}}}\rightarrow 0$ and we apply
Proposition~11(e)(ii).\vspace{-.5pc}
%\end{enumerate}

\item The fact that $U$ is a Dunford--Pettis operator follows from
Theorem~3.1 in \cite{10}, because $l_{1}$ has the Schur property.
Define $x_{n}^{*}\in l_{1}^{*}$ by $x_{n}^{*}(x) =\alpha
_{n}x_{n}$ and observe that, since $\alpha \in l_{\infty}$, then
the multiplication operator $ M\hbox{\rm :}\ l_{1}\rightarrow
l_{p}$ defined by $M(x) =(x_{n}^{*}(x))_{n\in \Bbb{N}}$ is
well-defined bounded linear. If $U$ is compact, then by
Proposition~11(c), $|\alpha_{n}|=\|x_{n}^{*}\|\rightarrow 0$.

\hskip 1pc Conversely, if $\alpha \in c_{0}$, then $M$ is compact
and hence by Proposition~11(a), $U$ is compact.

\item and (e). In our hypothesis, for any $1<p<\infty $ the
operator $ T\hbox{\rm :}\ C[0,1]\rightarrow l_{p}$ defined by
$T(f) = \big(\alpha_{n}\int_{0}^{1}f(s) r_{n}(s)
\hbox{d}s\big)_{n\in \Bbb{N}}$ is bounded linear (for $1<p<2$ by
the H\"{o}lder--Bessel inequality and $\alpha \in l_{r}$, while
for $2\leq p<\infty $ by $\|\|_{p}\leq \|\|_{2}$ and the Bessel
inequality).

\item Since $\alpha \in l_{r}$ for any $f\in C[0,1]$ we have the
inequality
\begin{equation*}
\hskip -1.25pc l_{r} \left(\alpha_{k}\int_{0}^{1}f(s) r_{k}(s)
\hbox{d}s| k\geq n \right) \leq l_{r}(\alpha_{k}| k\geq n) \|f\|.
\end{equation*}

From Exercise~6, p.~6 of \cite{4}, $T\hbox{\rm :}\
C[0,1]\rightarrow l_{r}$ is compact. Hence by the equivalences in
the proof of Proposition~11(d)(ii) $\big(\frac{1}{p}=\frac{
1}{2}+\frac{1}{r}\big)$, $U$ is compact.

\item Since $C([0,1],C[0,1])$ is isometric and isomorph with
$C([0,1]^{2})$ and $l_{p}$ is reflexive, $U$ is weakly compact and
hence is Dunford--Pettis (see Corollary~6, p.~154 of \cite{6}).

\hskip 1pc If $U$ is compact, then by Proposition~11(c),
$|\alpha_{n}|=\|x_{n}^{*}\|\rightarrow 0$.

\hskip 1pc For the converse, observe that for any $f\in C[0,1]$ we
have the following chain of inequalities:
\begin{align*}
\hskip -1.25pc l_{p} \left(\alpha_{k}\int_{0}^{1}f(s) r_{k}(s)
\hbox{d}s| k\geq n\right) &\leq (\sup\limits_{k\geq
n}|\alpha_{n}|) l_{p} \left(
\int_{0}^{1}f(s) r_{k}(s) \hbox{d}s| k\geq n \right)%\\
\end{align*}
\begin{align*}
\hskip -1.25pc &\leq (\sup\limits_{k\geq n}|\alpha_{n}|) l_{2}
\left(
\int_{0}^{1}f(s) r_{k}(s) \hbox{d}s| k\geq n \right)\\[.3pc]
\hskip -1.25pc &\leq (\sup\limits_{k\geq n}|\alpha_{n}|) \left(
\int_{0}^{1}
|f(s)|^{2}\hbox{d}s \right)^{1/2}\\[.3pc]
\hskip -1.25pc &\leq \|f\|(\sup\limits_{k\geq n}|\alpha_{n}|).
\end{align*}

Since $\alpha \in c_{0}$ and from Exercise~6, p.~6 of \cite{4}, it
follows that $T\hbox{\rm :}\ C[0,1]\rightarrow l_{p}$ is compact
and by Proposition~11(a), $U$ will be compact.

\item By the Hausdorff--Young inequality, the operator $T\hbox{\rm
:}\ L_{p^{*}}[0,1] \rightarrow l_{p}$ defined by $T(f)
=\big(\alpha_{n}\int_{0}^{1}f(s) r_{n}(s) \hbox{d}s\big)_{n\in
\Bbb{N}}$ is bounded linear.

\hskip 1pc Suppose $U$ is a Dunford--Pettis operator. Then by
Proposition~11(e)(i) $ T\hbox{\rm :}\ L_{p^{*}}[0,1] \rightarrow
c_{0}$ is a Dunford--Pettis operator and thus is compact, since
$L_{p^{*}}[0,1] $ is reflexive, hence
$|\alpha_{n}|=\|x_{n}^{*}\|\rightarrow 0$. The converse follows
from Proposition~11(e)(ii), since
$(L_{p^{*}}[0,1])^{*}=L_{p}[0,1]$ has type $\min (p,2) =2$ (see
Corollary~11.7, p.~219 of\break \cite{5}).

\item In our hypothesis the operator $T\hbox{\rm :}\ L_{1}[0,1]
\rightarrow l_{p}$ defined by $T(f) = \big(\int_{0}^{1}f(s)
g_{n}(s)$ $\hbox{d}s\big)_{n\in \Bbb{N}}$ is bounded linear and
weakly compact, and hence is a Dunford--Pettis \hbox{operator}
(see Lemma~4, p.~62 of \cite{6}) and Dunford--Pettis theorem
(p.~76). By Proposition~11(d)(i), $U$ is a Dunford--Pettis
operator. For the compactness we use Proposition~11(d)(ii) and the
representation of compact operators on $L_{1}(\mu)$, (see p.~68
of\break \cite{6}).\vspace{-1.5pc}
\end{enumerate}
\end{proof}

Now we analyze the case $l_{1}$. Since $l_{1}$ has the Schur
property (see Chapter~1 of \cite{5}), we study only the
compactness.

\begin{propo}$\left.\right.$\vspace{.5pc}

\noindent Let $X$ be a Banach space{\rm ,} $(x_{n}^{*})_{n\in
\Bbb{N}}\in w_{1}(X^{*})${\rm ,} $T\hbox{\rm :}\ X\rightarrow
l_{1}$ defined by $T(x) =(x_{n}^{*}(x))_{n\in \Bbb{N}}$ and let
$U\hbox{\rm :}\ C([0,1],X) \rightarrow l_{1}$ be the operator
defined by
\begin{equation*}
U(f) = \left(\int_{0}^{1}x_{n}^{*}f(t) r_{n}(t)\hbox{\rm d}t
\right)_{n\in \Bbb{N}}.
\end{equation*}

\begin{enumerate}
\renewcommand\labelenumi{\rm (\alph{enumi})}
\leftskip .15pc
\item
If the series $\sum_{n=1}^{\infty}x_{n}^{*}$ is unconditionally
norm convergent i.e. $T\hbox{\rm :}\ X\rightarrow l_{1}$ is a
compact operator{\rm ,} then $U$ is a compact operator.

\item If $U$ is a compact operator{\rm ,} then $T\hbox{\rm :}\
X\rightarrow l_{2}$ is a compact operator.
\end{enumerate}
\end{propo}

\begin{proof}$\left.\right.$

\begin{enumerate}
\renewcommand\labelenumi{\rm (\alph{enumi})}
\leftskip .15pc
\item It is analogous with we gave in Proposition~11(a).

\item If $U$ is compact, then, see Theorem~6 of \cite{1} and
Exercise~6, p.~6 of \cite{4}. For any $\varepsilon >0$ there is
$n_{\varepsilon }\in \Bbb{N}$ such that for any $n\geq
n_{\varepsilon}$ it follows that for any $f\in B(\Sigma,X) $ with
$\|f\|\leq 1$ we have
\begin{equation*}
\hskip -1.25pc \sum_{k=n_{\varepsilon}}^{n}
\left|\int_{0}^{1}x_{k}^{*}f(t)
r_{k}(t)\hbox{d}t\right|<\varepsilon.
\end{equation*}

Let $n\geq n_{\varepsilon}$. In particular, by the above
inequality, for any $ E\in \Sigma $ and $\|x\|\leq 1$ we have
\begin{equation*}
\hskip -1.25pc \left|\int_{E} \left(
\sum\limits_{k=n_{\varepsilon}}^{n} x_{k}^{*}(x) r_{k}(t) \right)
\hbox{d}t \right|<\varepsilon.
\end{equation*}

Then for any $\|x\|\leq 1$ we obtain
\begin{equation*}
\hskip -1.25pc \frac{1}{4}\int_{0}^{1} \left|
\sum\limits_{k=n_{\varepsilon}}^{n} x_{k}^{*}(x) r_{k}(t)
\right|\hbox{d}t\leq \sup\limits_{E\in \Sigma} \left| \int_{E}
\left( \sum\limits_{k= n_{\varepsilon}}^{n} x_{k}^{*}(x) r_{k}(t)
\right) \hbox{d}t \right|\leq \varepsilon
\end{equation*}
and, by Khinchin's inequality, $\frac{1}{4\sqrt{2}}
\big(\sum_{k=n_{\varepsilon}}^{n}|x_{k}^{*}(x)|^{2}\big)^{1/2}\leq
\varepsilon$, which by Exercise~6, p.~6 of \cite{4}, means that
$T\hbox{\rm :}\ X\rightarrow l_{2}$ is compact.\vspace{-1.8pc}
\end{enumerate}
\end{proof}

Now a concrete example.

\begin{exam}$\left.\right.$
{\rm
\begin{enumerate}
\renewcommand\labelenumi{\rm (\alph{enumi})}
\leftskip .15pc
\item Let $a=(a_{n})_{n\in \Bbb{N}}\in
l_{\infty}$ and let $U\hbox{\rm :}\ C([0,1],l_{1}) \rightarrow
l_{1}$ be the operator defined by
\begin{equation*}
\hskip -1.25pc U(f) = \left(a_{n}\int_{0}^{1}\langle
f(t),e_{n}\rangle r_{n}(t)\hbox{d}t \right)_{n\in \Bbb{N}}.
\end{equation*}
Then $U$ is a compact operator $\Leftrightarrow$ $a\in c_{0}$.

\item Let $g=(g_{n})_{n\in \Bbb{N}}\in L_{\infty}([0,1],l_{1})$
and $U\hbox{\rm :}\ C([0,1],L_{1}[0,1]) \rightarrow l_{1}$ be the
operator defined by
\begin{equation*}
\hskip -1.25pc U(f)= \left(\int_{0}^{1} \left(\int_{0}^{1}f(t) (s)
g_{n}(s) \hbox{d}s \right) r_{n}(t)\hbox{d}t \right)_{n\in
\Bbb{N}}.
\end{equation*}

%\begin{enumerate}
%\renewcommand\labelenumii{\rm (\roman{enumii})}
(i) If the function $g$ has an essentially relatively compact
range in $l_{1}$, then the operator $U$ is compact.

(ii) If $U$ is a compact operator, then the function $g\in
L_{\infty}([0,1],l_{1})$ has an essentially relatively compact
range in $l_{2}$.
%\end{enumerate}

\item Let $a=(a_{n})_{n\in \Bbb{N}}\in l_{\infty}$, $(E_{n})_{n\in
\Bbb{N}}\subseteq \Sigma$ be a sequence of pair-wise disjoint and
non-negligible Lebesgue sets. The operator $U\hbox{\rm :}\
C([0,1],L_{1}[0,1]) \rightarrow l_{1}$ defined by
\begin{equation*}
\hskip -1.25pc U(f) = \left( \int_{0}^{1}a_{n} \left(
\int_{E_{n}}f(t) (s) \hbox{d}s \right) r_{n}(t)\hbox{d}t
\right)_{n\in \Bbb{N}}
\end{equation*}
is a compact operator $\Leftrightarrow$ $a\in
c_{0}$.\vspace{-.5pc}
\end{enumerate}}
\end{exam}

\begin{proof}$\left.\right.$

\begin{enumerate}
\renewcommand\labelenumi{\rm (\alph{enumi})}
\leftskip .15pc
\item If $U$ is compact, then by Proposition~13(b) the
multiplication operator $M\hbox{\rm :}\ l_{1}\rightarrow l_{2}$
defined by $M(x_{n})_{n\in \Bbb{N}}=(a_{n}x_{n})_{n\in \Bbb{N}}$
is compact, which implies that $a\in c_{0}$. Conversely, if $a\in
c_{0}$, then the operator $M\hbox{\rm :}\ l_{1}\rightarrow l_{1}$
defined by $M(x_{n})_{n\in \Bbb{N}}=(a_{n}x_{n})_{n\in \Bbb{N}}$
is compact and thus, by Proposition~13(a), $U$ is compact.

\item
(i)~It follows from Remark~9(b) and Proposition~13(a). (ii)~If $U$
is a compact operator, then by Proposition~13(b), the operator
$T\hbox{\rm :}\ L_{1}[0,1] \rightarrow l_{2}$ defined by $T(f)
=\int_{0}^{1}f(t) g(t) $ is compact and using the representation
of compact operators on $L_{1}(\mu)$ in p.~68 of \cite{6}, the
statement follows.

\item By hypothesis $g=(a_{n}\chi_{E_{n}})_{n\in \Bbb{N}}\in
L_{\infty}([0,1],l_{1})$. If $U$ is compact, then by (b)(ii), $g$
has an essentially relatively compact range in $l_{2}$ i.e. there
is $A\in \Sigma$ with $\mu (A) =0$ such that for any $\varepsilon
>0$ there is $n_{\varepsilon}\in \Bbb{N}$ such that
$\sup_{t\notin A}\sum_{k=n_{\varepsilon}}^{\infty}
|a_{k}\chi_{E_{k}}(t)|^{2}<\varepsilon^{2}$. Let $n\geq
n_{\varepsilon}$. Since $E_{n}$ is a non-negligible Lebesgue,
there is $t\in E_{n}-A$ and thus $|a_{n}|^{2}<\varepsilon^{2}$
i.e. $a_{n}\rightarrow 0$.

\hskip 1pc The converse follows from inequality
$\sum_{k=n}^{\infty} |a_{n}\chi_{E_{k}}(t)|\leq \sup_{k\geq
n}|a_{k}|$ and (b)(i).\vspace{-2pc}
\end{enumerate}
\end{proof}

\end{document}